\documentclass[a4paper,11pt]{article}
\usepackage[top=30truemm,bottom=30truemm,left=30truemm,right=30truemm]{geometry}
\usepackage{amsmath}
\usepackage{amssymb}
\usepackage{float}
\usepackage{color}
\usepackage[dvipdfmx]{graphicx}
\usepackage{graphics}
\usepackage{graphicx}
\usepackage{tcolorbox}
\usepackage{latexsym}
\usepackage{ascmac}
\usepackage{framed}
\usepackage{fancybox}
\usepackage{theorem}
\usepackage{stmaryrd}
\usepackage{here}
\usepackage{mathtools}
\usepackage{appendix}

\mathtoolsset{showonlyrefs=true}

\theoremstyle{break}
\def\qed{\hfill$\Box$}

\newtheorem{defn}{Definition}
\newtheorem{lem}{Lemma}

\newtheorem{thm}{Theorem}

\newtheorem{rem}{Remark}

\makeatletter
\@addtoreset{equation}{section}

\makeatother

\makeatletter
\@addtoreset{figure}{section}
\def\thefigure{\thesection.\arabic{figure}}
\makeatother

\title{
Geometric structure of stationary problem for spatial 1D self-diffusion equation with logistic growth
}

\author{
Yu Ichida
\thanks{Department of Mathematical Sciences, School of Science, Kwansei Gakuin University, Gakuen Uegahara 1, Sanda, Hyogo 669-1330, Japan, {\tt ichi58yu@gmail.com} }
}

\begin{document}
\maketitle
\begin{abstract}
This paper considers the solution structure of non-trivial, non-constant stationary states of 1D spatial parabolic equations with nonlinear self-diffusion and logistic growth terms.
A two-dimensional ordinary differential equation satisfying the stationary problem is derived and all its dynamics, including to infinity, is revealed by the Poincar\'e-Lyapunov compactification, one of the compactifications of phase space.
The advantage of this method is that it can be used to classify all dynamical systems (especially connecting orbits) of a two-dimensional system including infinity.
Therefore, the classification results for the dynamical system including to infinity give the classification results for the non-constant stationary states obtained only from the structure of the original equations.
This argument allows us to observe a change in the classification of the non-constant stationary states by an explicit relation between the linear diffusion coefficient and the self-diffusion coefficient, combined with arguments about the symmetries and conserved quantities of the ODEs.
This means that changing the self-diffusion coefficient as a bifurcation parameter not only qualitatively changes the dynamical system from a big saddle homoclinic orbit of the ODEs to a heteroclinic orbit that connects the saddle equilibria, but also significantly changes the shape and the properties of the stationary states.
It explicitly shows the relationship between linear and self-diffusion, gives a characterization of non-trivial stationary states in terms of dynamical systems, and gives a deep insight into the influence of self-diffusion, one of the nonlinear diffusions.
\end{abstract}

{\bf Keywords:}
1D self-diffusion equation with logistic growth
stationary problem
Poincar\'e-type compactifcation, 
dynamics at infinity

\begin{center}
\end{center}

\section{Introduction}
\label{sec:SKTSS-int}
This paper, considers the following spatial one-dimensional parabolic equation with nonlinear self-diffusion and logistic growth terms:
\begin{equation}
u_{t}=\{(D+\alpha u) u\}_{xx}+\mu u(1-u), \quad t>0, \quad x\in \mathbb{R}, \quad u=u(t,x).
\label{eq:SKTSS-int1}
\end{equation}
Let $u_{t}=\partial u/\partial t$, $u_{xx}=\partial^{2}u/\partial x^{2}$ and $D$, $\alpha$ and $\mu$ are positive constants.

With $\Omega$ as the appropriate domain and with appropriate initial and boundary conditions in the equation
\begin{equation}
\begin{cases}
u_{t}=\Delta \{(D_{u}+\alpha_{1}u+\beta_{1}v) u\} + (r_{1}-a_{1}u-b_{1}v)u, \\
v_{t}=\Delta \{(D_{v}+\alpha_{2}u+\beta_{2}v) v\} + (r_{2}-a_{2}u-b_{2}v)v,
\end{cases}
\quad t>0, \quad x\in \Omega,
\label{eq:SKTSS-int2}
\end{equation}
this is called the SKT cross-diffusion equation.

The equation \eqref{eq:SKTSS-int1} is derived from
\begin{equation}
u_{t}=\Delta \{(D_{u}+\alpha_{1}u) u\} + (r_{1}-a_{1}u)u,
\label{eq:SKTSS-int3}
\end{equation}
which is obtained by considering $v=0$ and the whole domain in \eqref{eq:SKTSS-int2}.
In the SKT cross-diffusion equation\eqref{eq:SKTSS-int2}, $u=u(t,x)$ and $v=v(t,x)$ denote the density of the two species.
In addition, $D_{u}$ and $D_{v}$ are the diffusion coefficients of $u$ and $v$, respectively.
$r_{1}>0$ and $r_{2}>0$ are intrinsic growth rates, $a_{1}>0$ and $b_{2}>0$ are intra-specific competition rates, $a_{2}>0$ and $b_{1}>0$ are inter-specific competition rates.
The parameters $\alpha_{1}>0$ and $\beta_{2}>0$ indicate self-diffusion, while $\beta_{1}>0$ and $\alpha_{2}>0$ indicate cross-diffusion.

The equation \eqref{eq:SKTSS-int2} for the case of linear diffusion with $\alpha_{j}=\beta_{j}=0$ ($j=1,2$) is called the Lotka-Volterra reaction-diffusion system.
In this equation, it is known that stable stationary solutions describing biological segregation phenomena do not occur.
Therefore, the model \eqref{eq:SKTSS-int2} is a system of equations proposed by Shigesada-Kawasaki-Teramoto (\cite{SKT}) to describe biological segregation phenomena.
In addition, see \cite{MiKa, SKT2, SKT3} and references therein.
It is noteworthy that this proposal attempts to achieve its objective by adding a nonlinear diffusion term.
However, due to the complexity of the nonlinearity, it is difficult to directly answer the question of what happens to the solution profile given appropriate initial conditions, given the history of the SKT cross-diffusion equation.
Of course, the analysis of \eqref{eq:SKTSS-int2} is not easy even for stationary problems.
There are many approaches to stationary problems, such as the bifurcation theory, a singular limit and a singular perturbation technique and an elliptic approach.
For instance, see \cite{IzuKoba} and references therein.
Although many studies exist, the effects of cross-diffusion and self-diffusion on the shape of the solution are not successfully understood.
Therefore, it is important to investigate the effects of each and both cross-diffusion and self-diffusion to understand the behavior of the solution of the SKT cross-diffusion equation.
As a first step, we investigate the effect of nonlinearity only in the self-diffusion term, without considering the cross-diffusion term, using a compactification of the phase space, which is a completely different approach from the previous stationary problems.

Considering nonlinear diffusion with only a self-diffusion term implies the well-known structure of the equation in the following two senses.
The first is
\begin{equation}
u_{t}= Du_{xx}+\mu u(1-u),
\quad t>0, \quad x\in \mathbb{R}
\label{eq:SKTSS-int4}
\end{equation}
in the equation \eqref{eq:SKTSS-int1} with $\alpha\to0$, which is a Fisher-KPP type reaction-diffusion equation.
The second is
\begin{equation}
u_{t}=(\alpha u^{2})_{xx}+\mu u(1-u), 
\quad t>0, \quad x\in \mathbb{R}
\label{eq:SKTSS-int5}
\end{equation}
with $D\to0$ in the equation \eqref{eq:SKTSS-int1}, which is a degenerate reaction-diffusion equation of the Porous-Fisher-KPP type with nonlinear diffusion.
The equation \eqref{eq:SKTSS-int1} is a reduced equation for SKT cross diffusion such that in addition to including both linear and nonlinear diffusion, it also includes a typical reaction-diffusion equation and a nonlinear diffusion equation of the porous medium type.
This means that the question naturally arises as to how the existence, shape, and asymptotic behavior of the solution change with the influence of $D$ and $\alpha$.
The approach taken in this paper is to focus on stationary problems and to classify information about the existence and profiles of the functions that satisfy them from the structure of the equations alone.
We expect that the enumeration of this information will lead to an understanding of how the solution evolves as an equation depending on the initial conditions.

This paper aims to classify information on the existence and shape of functions to be satisfied by stationary problems in \eqref{eq:SKTSS-int1}.
The author has previously clarified the structure of special solutions of partial differential equations in \cite{QTW, cDNPE, pDNPE, BIRD, IM-DNPE} by using Poincar\'e-type compactification, a dynamical systems theory and a geometric approach, focusing only on the structure of the equations.
The Poincar\'e-type compactification plays a central role in this paper, see \cite{FAL, QTW, BIRD, Matsue1, Matsue2} or Appendix A for details.
This method is one of the compactifications of the original phase space and is the embedding of $\mathbb{R}^{n}$ in $\mathbb{R}^{n+1}$ in the unit upper hemisphere.
This paper also considers a two-dimensional system of ordinary differential equations (hereafter ODEs) corresponding to the case $n=2$.
In the case that $n$ is $n=2$, the part of the phase space corresponding to infinity is $\mathbb{S}^{1}$. 
In order to capture this relatively easily in the sense of a transformation of differential equations, the system is divided into several parts covering $\mathbb{S}^{1}$ and each is considered by projecting it into local coordinates.
By integrating the dynamical system in these local coordinates, we can understand the behavior at infinity.
This method has been used, for instance, in the analysis of the Li\'enard equation (see \cite{FAL} and references therein) and in the dynamical system theory reconstruction of the blow-up solution of ODEs (see \cite{Matsue1, Matsue2}).

The stationary problem for the equation\eqref{eq:SKTSS-int1} is
\begin{equation}
0=[(D+2\alpha u)u_{x}]_{x}+\mu u(1-u), \quad x\in \mathbb{R}, \quad u=u(x),
\label{eq:SKTSS-int6}
\end{equation}
and the following transformation is introduced to make the ODE \eqref{eq:SKTSS-int6} easier to handle:
\begin{equation}
v(x) := (D+2\alpha u)u_{x}.
\label{eq:SKTSS-int7}
\end{equation}
Then, the following two-dimensional ODEs are obtained:
\begin{equation}
\begin{cases}
u_{x}= (D+2\alpha u)^{-1}v, \\
v_{x}= -\mu u(1- u)
\end{cases}
\label{eq:SKTSS-int8}
\end{equation}

In this paper, we study the dynamical systems of \eqref{eq:SKTSS-int8} including to infinity by the Poincar\'e-type compactification, and the classification of all dynamical systems including to infinity corresponds to the classification of functions satisfying stationary problems.
One of the conclusions of this paper is that the structure of the nonconstant stationary states can be classified into three types for $D<2\alpha$, $D=2\alpha$, and $D>2\alpha$.
This change in classification can then be characterized by a bifurcation.
For $D=2\alpha$, the orbits are heteroclinic between different saddles, while for $D\neq 2\alpha$, they are ``big'' saddle homoclinic orbits (for instance, see \cite{Kuz} and references therein) in which the orbit starts from a saddle and returns to itself.
This fact is derived from the symmetry and the conserved quantities of the equation \eqref{eq:SKTSS-pd2} that is tranformed by time-scale desingularization.
This means that it occurs the bifurcation by the relationship between the self-diffusion coefficient $\alpha$ and the linear diffusion coefficient $D$, and that this bifurcation significantly changes the characteristics of the stationary states.
The change in the self-diffusion coefficient $\alpha$ is discussed in detail in Section \ref{sec:SKTSS-di}.

The organization of the paper is organized as follows.
In the next section, we present the main results of this paper for the stationary problem \eqref{eq:SKTSS-int6} of \eqref{eq:SKTSS-int1}.
Section \ref{sec:SKTSS-pd} gives the classification of the dynamical system at $\mathbb{R}^{2}\, \cup\, \{\|(\phi, \psi)\|=+\infty\}$, which corresponds to the dynamics of \eqref{eq:SKTSS-int8} including to infinity.
This dynamical system including to infinity is also called the Poincar\'e-Lyapunov disk, which is obtained by the Poincar\'e-Lyapunov compactification.
Section \ref{sec:SKTSS-pro} gives a proof of the main results.
Finally, we give some discussion of the results on the classification of functions to be satisfied for stationary problems with nonlinear diffusion coefficients, including the equations considered in this paper in Section \ref{sec:SKTSS-di}.

\section{Main results}
\label{sec:SKTSS-mr}
The main results of this paper are described.
There exist constant stationary solutions of the equation \eqref{eq:SKTSS-int6} such that $u\equiv 0$ and $u\equiv 1$.
This paper focuses on nontrivial nonconstant stationary states such that $u(x)\ge -(2\alpha)^{-1}D$ is satisfied and gives classification results for them.
In the following, the symbol $f(x)\sim g(x)$ as $x\to a$ implies that
\[
\lim_{x\to a}\left| \dfrac{f(x)}{g(x)} \right|=1.
\]
For the classification of functions satisfying stationary problems, except for periodic ones, we classify the types satisfied by the stationary states in terms of their behavior in $x \to +\infty$(or $x\to x_{+}-0$ as $|x_{+}|<+\infty$) and in $x \to -\infty$ (or $x\to x_{-}+0$ as $|x_{-}|<+\infty$) by introducing several results, and in combination give a classification on the existence of stationary states and their characterization.

\begin{defn}
\label{def:SKTSS-mr1}
Assume that $D, \alpha, \mu$ are positive constants.
Each function $u(x)$ such that it satisfies \eqref{eq:SKTSS-int6} is defined as follows according to its behavior at $x\to +\infty$:
\begin{enumerate}
\item[(i)]
Define a function of type $*-1$ such that 
\[
\lim_{x \to +\infty} u(x) =1, \quad  \lim_{x \to +\infty} u_{x}(x)=0
\]
is satisfied as $x\to +\infty$.
It corresponds to the orbit in \eqref{eq:SKTSS-int8} at $(u, v)=(1, 0)$.
\item[(ii)]
Define an unbounded function of type $*-\infty$ such that 
\[
\lim_{x \to +\infty} u(x) = \lim_{x \to +\infty} u_{x}(x)=+\infty
\]
is satisfied as $x\to +\infty$.
It corresponds to the orbit in \eqref{eq:SKTSS-int8} at $(u, v)=(+\infty, +\infty)$.
\end{enumerate}
Similarly, each function $u(x)$ such that it satisfies \eqref{eq:SKTSS-int6} is defined as follows according to its behavior at $x\to -\infty$:
\begin{enumerate}
\item[(iii)]
Define a function of type $1-*$ such that 
\[
\lim_{x \to -\infty} u(x) =1, \quad  \lim_{x \to -\infty} u_{x}(x)=0
\]
is satisfied as $x\to -\infty$.
It corresponds to the orbit in \eqref{eq:SKTSS-int8} at $(u, v)=(1, 0)$.
\item[(iv)]
Define an unbounded function of type $\infty-*$ such that 
\[
\lim_{x \to -\infty} u(x) =+\infty, \quad  \lim_{x \to -\infty} u_{x}(x)=-\infty
\]
is satisfied as $x\to -\infty$.
It corresponds to the orbit in \eqref{eq:SKTSS-int8} at $(u, v)=(+\infty, -\infty)$.
\end{enumerate}
\end{defn}

Among the functions that satisfy the stationary problem in this paper, there are functions that are satisfied in a finite or semi-infinite interval in \eqref{eq:SKTSS-int6} and have singularity at the endpoints.
However, we consider the whole domain.
The treatment of the solution is not easy.
Therefore, the following definitions are adopted to handle these functions.
These definitions, which are also given in the traveling wave case in \cite{cDNPE, pDNPE, IM-DNPE}, are reproduced for the convenience of the reader.

\begin{defn}
\label{def:SKTSS-mr2}
We say that a function $u(x)$ is a quasi-stationary state of \eqref{eq:SKTSS-int1} if the function $u(x)$ is a solution of \eqref{eq:SKTSS-int6} on a finite interval or semi-infinite interval.
\end{defn}

\begin{defn}
\label{def:SKTSS-mr3}
Assume that $D, \alpha, \mu$ are positive constants.
\begin{enumerate}
\item[(v)]
Define a function of type $qs_{-}-*$ such that the function $u(x)$ is a quasi-stationary state and it satisfies 
\[
\lim_{x \to x_{-} + 0} u(x) = -\dfrac{D}{2\alpha}, \quad
\lim_{x \to x_{-} +0} u_{x}(x) =  C
\]
with $|x_{-}|<\infty$ and a positive constant $C$.
It corresponds to the orbit in \eqref{eq:SKTSS-int8} at $(u, v)=(-(2\alpha)^{-1}D, 0)$. 
\item[(vi)]
Define a function of type $*-qs_{+}$ such that the function $u(x)$ is a quasi-stationary state and it satisfies 
\[
\lim_{x \to x_{+} - 0} u(x) = -\dfrac{D}{2\alpha}, \quad
\lim_{x \to x_{+} -0} u_{x}(x) =  -C
\]
with $|x_{+}|<\infty$ and a positive constant $C$.
It corresponds to the orbit in \eqref{eq:SKTSS-int8} at $(u, v)=(-(2\alpha)^{-1}D, 0)$. 
\end{enumerate}
\end{defn}

\begin{thm}
\label{th:SKTSS-mr1}
Assume that $D, \alpha, \mu$ are positive constants.
For any given $D, \alpha, \mu$, the equation \eqref{eq:SKTSS-int6} has four types of functions $u(x)$.
\begin{enumerate}
\item[(I)] 
There exists an unbounded stationary solution of type $1-\infty$.
In addition, the asymptotic behavior of $u(x)$ for $x \to +\infty$ is
\begin{equation}
u(x) \sim B_{1}e^{\frac{\sqrt{\alpha\mu}}{2\alpha}x}
\quad {\rm{as}} \quad x \to +\infty
\label{eq:SKTSS-mr1}
\end{equation}
and the asymptotic behavior of $u(x)$ for $x \to -\infty$ is
\begin{equation}
u(x) \sim 1+B_{2}e^{\frac{\omega_{+}}{2\alpha+D}x}
\quad {\rm{as}} \quad x \to -\infty,
\label{eq:SKTSS-mr2}
\end{equation}
where $B_{1,2}$ are positive constants and $\omega_{+}= \sqrt{\mu(2\alpha+D)}$ holds.
\item[(II)] 
There exists a family of unbounded stationary solutions of type $\infty - \infty$.
There exists a constant $x_{0}\in (-\infty, +\infty)$ such that the following holds:
$u_{x}(x)<0$ for $x\in (-\infty, x_{0})$, $u_{x}(x_{0})=0$ and $u_{x}(x)>0$ for $x\in (x_{0}, +\infty)$.
In addition, the asymptotic behavior of $u(x)$ for $x \to +\infty$ is \eqref{eq:SKTSS-mr4} and the asymptotic behavior of $u(x)$ for $x \to -\infty$ is
\begin{equation}
u(x) \sim B_{3}e^{-\frac{\sqrt{\alpha\mu}}{2\alpha}x}
\quad {\rm{as}} \quad x \to -\infty,
\label{eq:SKTSS-mr3}
\end{equation}
where $B_{3}$ is a positive constant.
\item[(III)] 
There exists an unbounded stationary solution of type $\infty-1$.
In addition, the asymptotic behavior of $u(x)$ for $x \to -\infty$ is \eqref{eq:SKTSS-mr6} and the asymptotic behavior of $u(x)$ for $x \to +\infty$ is
\begin{equation}
u(x) \sim 1+B_{4}e^{\frac{\omega_{-}}{2\alpha+D}x}
\quad {\rm{as}} \quad x \to +\infty,
\label{eq:SKTSS-mr4}
\end{equation}
where $B_{4}$ is a positive constant and $\omega_{-}= -\sqrt{\mu(2\alpha+D)}$ holds.
\item[(IV)]
There exists a family of stationary solutions with periodicity such that it corresponds to the family of periodic orbits in \eqref{eq:SKTSS-int8}.
\end{enumerate}
\end{thm}

\begin{thm}
\label{th:SKTSS-mr2}
Assume that $0<D, \alpha, \mu \in\mathbb{R}$.
If  $D<2\alpha$, in addition to the functions obtained in Theorem \ref{th:SKTSS-mr1},  the equation \eqref{eq:SKTSS-int6} has a function of type $qs_{-}-qs_{+}$.
It holds $u_{x}(x)>0$ for $x \in (-\infty, x^{*})$, $u_{x}(x^{*})=0$ and $u_{x}(x)<0$ for $x\in (x^{*}, +\infty)$ and $0<u(x^{*})<1$.
In addition, the asymptotic behavior of $u(x)$ for $x \to x_{-}+0$ is
\begin{equation}
u(x) \sim A_{1}(x-x_{-})-\dfrac{D}{2\alpha}
\quad {\rm{as}} \quad x \to x_{-}+0
\label{eq:SKTSS-mr5} 
\end{equation}
and the asymptotic behavior of $u(x)$ for $x \to x_{+}-0$ is
\begin{equation}
u(x) \sim A_{2}(x_{+}-x)-\dfrac{D}{2\alpha}
\quad {\rm{as}} \quad x \to x_{+}-0,
\label{eq:SKTSS-mr6} 
\end{equation}
where $A_{1,2}>0$ are constants.
\end{thm}

\begin{thm}
\label{th:SKTSS-mr3}
Assume that $0<D, \alpha, \mu \in\mathbb{R}$.
If  $D=2\alpha$, in addition to the functions obtained in Theorem \ref{th:SKTSS-mr1}, the equation \eqref{eq:SKTSS-int6} has two types as follows:
\begin{enumerate}
\item[(a)]
The equation \eqref{eq:SKTSS-int6} has a function of type $1-qs_{+}$.
In addition, the asymptotic behavior of $u(x)$ for $x\to -\infty$ is
\begin{equation}
u(x) \sim 1-A_{3}e^{\frac{\omega_{+}}{2\alpha+D}x}
\quad {\rm{as}} \quad x \to -\infty
\label{eq:SKTSS-mr7}
\end{equation}
with $A_{3}>0$ is constant and the asymptotic behavior of $u(x)$ for $x \to x_{+}-0$ is \eqref{eq:SKTSS-mr6}.
\item[(b)]
The equation \eqref{eq:SKTSS-int6} has a function of type $qs_{-}-1$.
In addition, the asymptotic behavior of $u(x)$ for $x \to x_{-}+0$ is \eqref{eq:SKTSS-mr5} and the asymptotic behavior of $u(x)$ for $x\to +\infty$ is
\begin{equation}
u(x) \sim 1-A_{4}e^{\frac{\omega_{-}}{2\alpha+D}x}
\quad {\rm{as}} \quad x \to +\infty
\label{eq:SKTSS-mr8}
\end{equation}
where $A_{4}>0$ is constant.
\end{enumerate}
\end{thm}

\begin{thm}
\label{th:SKTSS-mr4}
Assume that $0<D, \alpha, \mu \in\mathbb{R}$.
If  $D>2\alpha$, in addition to the functions obtained in Theorem \ref{th:SKTSS-mr1}, the equation \eqref{eq:SKTSS-int6} has three types as follows:
\begin{enumerate}
\item[(A)]
The equation \eqref{eq:SKTSS-int6} has a stationary solution of type $1-1$.
It holds $u_{x}(x)<0$ for $x \in (-\infty, x_{*})$, $u_{x}(x_{*})=0$ and $u_{x}(x)>0$ for $x\in (x_{*}, +\infty)$.
The asymptotic behavior of $u(x)$ for $x \to -\infty$ is \eqref{eq:SKTSS-mr7} and the asymptotic behavior of $u(x)$ for $x \to +\infty$ is \eqref{eq:SKTSS-mr8}.
\item[(B)]
The equation \eqref{eq:SKTSS-int6} has an unbounded sign-changing function of type $\infty-qs_{+}$.
In addition, the asymptotic behavior of $u(x)$ for $x \to -\infty$ is \eqref{eq:SKTSS-mr2} and the asymptotic behavior of $u(x)$ for $x \to x_{+}-0$ is \eqref{eq:SKTSS-mr6}.
\item[(C)]
The equation \eqref{eq:SKTSS-int6} has an unbounded sign-changing function of  type $qs_{-}-\infty$.
In addition, the asymptotic behavior of $u(x)$ for $x \to x_{-}-0$ is \eqref{eq:SKTSS-mr5} and the asymptotic behavior of $u(x)$ for $x \to +\infty$ is \eqref{eq:SKTSS-mr4}.
\end{enumerate}
\end{thm}

\section{Dynamics on the Poincar\'e-Lyapunov disk}
\label{sec:SKTSS-pd}
In this section, by using the Poincar\'e-Lyapunov compactification, we study the dynamics on the Poincar\'e-Lyapunov disk $\mathbb{R}^{2} \cup \{ (u, v) \mid \|(u, v)\|=+\infty \}$.

\subsection{Dynamics near finite equilibria}
\label{sub:SKTSS-pd1}
We desingularize $D+2\alpha u=0$ by the time-rescale desingularization:
\begin{equation}
d\tau/dx=(D+2\alpha u)^{-1}
\label{eq:SKTSS-pd1}
\end{equation}
as in \cite{cDNPE, pDNPE, sDNPE, IM-DNPE}.
Since we are considering on $\{u(x) \ge -(2\alpha)^{-1}D\}$, the direction of the time does not change via the desingularization \eqref{eq:SKTSS-pd1} in this region.
Then we have
\begin{equation}
\begin{cases}
u' = v, \\
v' = -\mu u(1-u)(D+2\alpha u), 
\end{cases}
\quad \left(\,\, '=\dfrac{d}{d\tau} \right).
\label{eq:SKTSS-pd2}
\end{equation}

\begin{rem}
\label{rem:SKTSS-pd1}
It should be noted that the time scale desingularization \eqref{eq:SKTSS-pd1} is simply multiplying the vector field by $D+2\alpha u$.
Then, except for the singularity $\{D+2\alpha u=0\}$, the integral curves of the system (vector field) remain the same but are parameterized differently in $\{u\ge -(2\alpha)^{-1}D\}$. 
We refer to section 7.7 of \cite{CK} and references therein for the analytical treatments of desingularization with the time rescaling. 
In what follows, we use a similar time rescaling (re-parameterization of the solution curves) repeatedly to desingularize the vector fields.
\end{rem}

\begin{rem}
\label{rem:SKTSS-pd2}
The equation \eqref{eq:SKTSS-pd2} coincides with the two-dimensional ODEs \eqref{eq:SKTSS-int8}, which appear as stationary problems for the equation
\begin{equation}
u_{t}= u_{xx}+ \tilde{\mu} u(1-u)(1+ku),
\quad t>0, \quad x\in \mathbb{R}, \quad u=u(t,x)
\label{eq:SKTSS-int9}
\end{equation}
in \cite{HaRo}.
However, $\tilde{\mu}=\mu D$, $k=2\alpha/D$.
In the sense of a nonconstant stationary solution, the nonconstant stationary solution of \eqref{eq:SKTSS-pd1} has a more complex solution structure since we use of the time scale transformation \eqref{eq:SKTSS-pd1}.
The classification of functions satisfying the stationary problem for the equation \eqref{eq:SKTSS-int9} is described in Appendix B.
\end{rem}

The system \eqref{eq:SKTSS-pd2} has the following equilibria:
\[
E_{0}: (u, v)=(0,0), \quad 
E_{1}: (u, v)=(1,0), \quad
E_{2}: (u, v)=(-(2\alpha)^{-1}D, 0).
\]
The Jacobian matrices $J_{i}$ ($i=0,1,2$) of the vector field \eqref{eq:SKTSS-pd2} at $E_{i}$ ($i=0, 1, 2$) are
\begin{align*}
& J_{0}=\left( \begin{array}{cc} 
0 & 1 \\ -\mu D & 0
\end{array} \right), \\
&J_{1}= \left( \begin{array}{cc} 
0 & 1 \\ 2\alpha\mu+\mu D & 0
\end{array} \right), \\
&J_{2}=\left( \begin{array}{cc} 
0 & 1 \\ (2\alpha)^{-1}\mu D[2\alpha+D] & 0
\end{array} \right).
\end{align*}
The eigenvalues and eigenvectors of each linearized matrix are organized as follows:
\begin{itemize}
\item 
The eigenvalue of $J_{0}$ is $\pm \sqrt{\mu D}i$.
Note that $i$ denotes imaginary units.
\item
The eigenvalue of $J_{1}$ is
\[
\omega_{\pm}=\pm \sqrt{\mu(2\alpha+D)})
\]
and the corresponding eigenvectors are $(1, \omega_{-})^{T}$ and $(1, \omega_{+})^{T}$, respectively
Here, $T$ is the symbol for transpose.
Therefore, $E_{1}$ is a saddle.
\item
The eigenvalue of $J_{2}$ is
\[
\Lambda_{\pm}=\pm \sqrt{\dfrac{\mu D^{2}+2\alpha \mu D}{2\alpha}}
\]
and the corresponding eigenvectors are $(1, \Lambda_{-})^{T}$ and $(1, \Lambda_{+})^{T}$, respectively.
Therefore, $E_{2}$ is a saddle.
\end{itemize}

\subsection{Asymptotically quasi-homogeneous vector field}
\label{sub:SKTSS-pd2}
We study the type and order of the vector field \eqref{eq:SKTSS-pd2} in preparation for studying the dynamics on $\mathbb{R}^{2} \cup \{(u, v) \mid \|(u, v)\|=+\infty\}$ by using the Poincar\'e-Lyapunov compactification.
This is an important preparation for studying the infinity dynamics of \eqref{eq:SKTSS-pd2}.
Note that proper study of these means choosing a compactification that reflects the original vector field information. 
In this paper, the Poincar\'e-Lyapunov compactification is used to properly extract the structure at infinity. 
See Appendix A for the definition in asymptotically quasi-homogeneous vector fields
See also \cite{FAL, Matsue1, Matsue2} and references therein.
Let $f=(f_{1}(u, v), f_{2}(u, v))$ be 
\[
f_{1}(u, v)= v,
\quad
f_{2}(u, v)= -\mu u(1-u)(D+2\alpha u).
\]
Then we have the following observation (see Appendix A and \cite{Matsue1, Matsue2}) for more details).

\begin{lem}
\label{lem:SKTSS-pd1}
Assume that $0<\mu, \alpha, D\in \mathbb{R}$.
Then, the vector field $f$ is asymptotically quasi-homogeneous of type $(1, 2)$ and order $2$ at infinity.
\end{lem}
{\bf{Proof.}}
Let a type be $(\alpha_{1}, \alpha_{2})$ and $R\in \mathbb{R}$.
For all $(u, v)\in \mathbb{R}^{2}$, the following holds:
\begin{align*}
f_{1}(R^{\alpha_{1}}u, R^{\alpha_{2}}v) &= R^{k+\alpha_{1}}f_{1}(u, v), \\
f_{2}(R^{\alpha_{1}}u, R^{\alpha_{2}}v) &= R^{k+\alpha_{2}}f_{2}(u, v).
\end{align*}
Using \eqref{eq:SKTSS-pd2}, the left-hand sides above are calculated as follows:
\begin{align*}
f_{1}(R^{\alpha_{1}}u, R^{\alpha_{2}}v) &=R^{\alpha_{2}}v,
\\
f_{2}(R^{\alpha_{1}}u, R^{\alpha_{2}}v) &= -\mu D R^{\alpha_{1}}u-\mu(2\alpha-D)R^{2\alpha_{1}}u^{2}+2\alpha\mu R^{3\alpha_{1}}u^{3}.
\end{align*}
By comparing the order parts, we get
\begin{equation}
\begin{cases}
\alpha_{2} = k+\alpha_{1},
\\
\alpha_{1}=k+\alpha_{2},
\\
2\alpha_{1}=k+\alpha_{2},
\\
3\alpha_{1}=k+\alpha_{2}
\end{cases}
\label{eq:SKTSS-pd3}
\end{equation}
Since the first and fourth equations in \eqref{eq:SKTSS-pd3} correspond to the maximum order in \eqref{eq:SKTSS-pd2}, $(\alpha_{1}, \alpha_{2})=(1,2)$ and $k = 1$ are obtained from them. 
These then satisfy the second equation in \eqref{eq:SKTSS-pd3}. 
Furthermore, we set 
\[
(f_{\alpha,k})_{1}= v,
\quad
(f_{\alpha,k})_{2}= 2\alpha\mu u^{3}
\]
and we then have 
\begin{align*}
&\lim_{R\to +\infty}R^{-(k+\alpha_{1})} \left\{ f_{1}(R^{\alpha_{1}}u, R^{\alpha_{2}}v) -R^{k+\alpha_{1}}(f_{\alpha,k})_{1}(u, v) \right\} =0,
\\
&\lim_{R\to +\infty}R^{-(k+\alpha_{2})} \left\{ f_{2}(R^{\alpha_{1}}u, R^{\alpha_{2}}v) -R^{k+\alpha_{2}}(f_{\alpha,k})_{2}(u, v) \right\} =0.
\end{align*}
From the above results, we can see that the vector field $f$ in \eqref{eq:SKTSS-pd2} is asymptotically quasi-homogeneous of type $(\alpha_{1}, \alpha_{2})=(1, 2)$ and order $k+1=2$ at infinity.
Thus, the proof of Lemma \ref{lem:SKTSS-pd1} is completed.
\qed
\\

We will consider the dynamics of this equation on the charts $\overline{U}_{j}$ and $\overline{V}_{j}$ ($j =1,2$). 
See Appendix A for the definitions and geometric images of $\overline{U}_{j}$ and $\overline{V}_{j}$, which is described similarly in the previous studies \cite{cDNPE, pDNPE, sDNPE, IM-DNPE}. 
See also \cite{FAL, Matsue1, Matsue2}.

\subsection{Dynamics on the local chart $\overline{U}_{1}$}
\label{sub:SKTSS-pd3}
To obtain the dynamics on the chart $\overline{U}_{1}$, we introduce the coordinates $(\lambda_{1}, \lambda_{2})$ given by 
\[
u=1/\lambda_{1}, \quad v=\lambda_{2}/\lambda_{1}^{2}.
\]
Here, note that the exponents of $\lambda_{1}$ are derived from the type found in Lemma \ref{lem:SKTSS-pd1}.
We then have 
\begin{equation}
\begin{cases}
\lambda_{1}'= -\lambda_{2},
\\
\lambda_{2}'=-\mu D \lambda_{1} -\mu(2\alpha-D)+2\alpha\mu \lambda_{1}^{-1}-2\lambda_{1}^{-1}\lambda_{2}^{2},
\end{cases}
\label{eq:-SKTSS-pd4}
\end{equation}
By using the time-scale desingularization $ds/d\tau=\lambda_{1}^{-1}$, we can obtain 
\begin{equation}
\begin{cases}
d\lambda_{1}/ds=-\lambda_{1} \lambda_{2},
\\
d\lambda_{2}/ds=-\mu D \lambda_{1}^{2} -\mu(2\alpha-D)\lambda_{1}+2\alpha\mu-2\lambda_{2}^{2}.
\end{cases}
\label{eq:SKTSS-pd5}
\end{equation}
The equilibria of the system \eqref{eq:SKTSS-pd5} on $\{\lambda_{1}=0\}$ are
\[
E_{3}: (\lambda_{1}, \lambda_{2})=(0,-\sqrt{\alpha\mu}), \quad
E_{4}: (\lambda_{1}, \lambda_{2})=(0, \sqrt{\alpha\mu}).
\]
The Jacobian matrices of the vector field \eqref{eq:SKTSS-pd5} at $E_{3, 4}$ are
\begin{align*}
& E_{3}: \left( \begin{array}{cc} 
\sqrt{\alpha\mu} & 0 \\ -\mu(2\alpha-D) & 4\sqrt{\alpha\mu}
\end{array} \right), \\
&E_{4}: \left( \begin{array}{cc} 
-\sqrt{\alpha\mu} & 0 \\ -\mu(2\alpha-D) & -4\sqrt{\alpha\mu}
\end{array} \right).
\end{align*}
Therefore, $E_{3}$ is an unstable equilibrium and $E_{4}$ is an asymptotically stable equilibrium.

\subsection{Dynamics on the local chart $\overline{V}_{1}$}
\label{sub:SKTSS-pd4}
The change of coordinates
\[
u=-1/\lambda_{1}, \quad
v=-\lambda_{2}/\lambda_{1}^{2}
\]
and the time-rescaling $ds/d\tau=\lambda_{1}^{-1}$ give the projection dynamics of \eqref{eq:SKTSS-pd2} on the chart $\overline{V}_{1}$.
Then we obtain the following system: 
\begin{equation}
\begin{cases}
d\lambda_{1}/ds=-\lambda_{1} \lambda_{2},
\\
d\lambda_{2}/ds= -\mu D \lambda_{1}^{2} +\mu(2\alpha-D)\lambda_{1}+2\alpha\mu-2\lambda_{2}^{2}.
\end{cases}
\label{eq:SKTSS-pd6}
\end{equation}
The equilibria of the system \eqref{eq:SKTSS-pd6} on  $\{\lambda_{1}=0\}$ are
\[
E_{5}: (\lambda_{1}, \lambda_{2})=(0,-\sqrt{\alpha\mu}), \quad
E_{6}: (\lambda_{1}, \lambda_{2})=(0, \sqrt{\alpha\mu})
\]
The Jacobian matrices of the vector field \eqref{eq:SKTSS-pd6} at $E_{5, 6}$ are
\begin{align*}
& E_{5}: \left( \begin{array}{cc} 
\sqrt{\alpha\mu} & 0 \\ \mu(2\alpha-D) & 4\sqrt{\alpha\mu}
\end{array} \right), \\
&E_{6}: \left( \begin{array}{cc} 
-\sqrt{\alpha\mu} & 0 \\ \mu(2\alpha-D) & -4\sqrt{\alpha\mu}
\end{array} \right)
\end{align*}
Therefore, $E_{5}$ is an unstable equilibrium and $E_{6}$ is an asymptotically stable equilibrium.

\subsection{Dynamics on the local chart $\overline{U}_{2}$}
\label{sub:SKTSS-pd5}
The change of coordinates $u=\lambda_{2}/\lambda_{1}$ and $v=1/\lambda_{1}^{2}$ and time-rescaling $ds /d\tau = \lambda_{1}^{-1}$ to study the dynamics on the local chart $\overline{U}_{2}$ yields
\begin{equation}
\begin{cases}
d\lambda_{1}/ds= 2^{-1}\mu D\lambda_{1}^{3}\lambda_{2}+2^{-1}\mu(2\alpha-D)\lambda_{1}^{2}\lambda_{2}^{2}-\alpha\mu\lambda_{1} \lambda_{2}^{3},
\\
d\lambda_{2}/ds=1+2^{-1}\mu D\lambda_{1}^{2}\lambda_{2}^{2}+2^{-1}\mu(2\alpha-D)\lambda_{1} \lambda_{2}^{3}-\alpha\mu \lambda_{2}^{4}.
\end{cases}
\label{eq:SKTSS-pd7}
\end{equation}
The equilibria of the system \eqref{eq:SKTSS-pd7} on $\{\lambda_{1}=0\}$ are
\[
E_{7}: (\lambda_{1}, \lambda_{2})=(0,-M), \quad
E_{8}: (\lambda_{1}, \lambda_{2})=(0, M), \quad
M= \sqrt{(\alpha\mu)^{-1}}.
\]
The Jacobian matrices of the vector field \eqref{eq:SKTSS-pd7} at $E_{7, 8}$ are
\begin{align*}
& E_{7}: \left( \begin{array}{cc} 
\alpha\mu M^{3} & 0 \\ -2^{-1}\mu(2\alpha-D)M^{3} & 4\alpha\mu M^{3}
\end{array} \right), \\
&E_{8}: \left( \begin{array}{cc} 
-\alpha\mu M^{3} & 0 \\ 2^{-1}\mu(2\alpha-D)M^{3} & -4\alpha\mu M^{3}
\end{array} \right)
\end{align*}
Therefore, $E_{7}$ is an unstable equilibrium and $E_{8}$ is an asymptotically stable equilibrium.

\subsection{Dynamics on the local chart $\overline{V}_{2}$}
\label{sub:SKTSS-pd6}
The change of coordinates $u=-\lambda_{2}/\lambda_{1}$ and $v=-1/\lambda_{1}^{2}$ and time-rescaling $ds /d\tau = \lambda_{1}^{-1}$ to study the dynamics on the local chart $\overline{V}_{2}$ yields
\begin{equation}
\begin{cases}
d\lambda_{1}/ds= 2^{-1}\mu D\lambda_{1}^{3}\lambda_{2} - 2^{-1}\mu(2\alpha-D)\lambda_{1}^{2}\lambda_{2}^{2}-\alpha\mu\lambda_{1} \lambda_{2}^{3},
\\
d\lambda_{2}/ds=1 - 2^{-1}\mu D\lambda_{1}^{2}\lambda_{2}^{2}-2^{-1}\mu(2\alpha-D)\lambda_{1} \lambda_{2}^{3}-\alpha\mu \lambda_{2}^{4}.
\end{cases}
\label{eq:SKTSS-pd8}
\end{equation}
The equilibria of the system \eqref{eq:SKTSS-pd8} on $\{\lambda_{1}=0\}$ are
\[
E_{9}: (\lambda_{1}, \lambda_{2})=(0,-M), \quad
E_{10}: (\lambda_{1}, \lambda_{2})=(0, M), \quad
\]
The Jacobian matrices of the vector field \eqref{eq:SKTSS-pd8} at $E_{9, 10}$ are
\begin{align*}
& E_{9}: \left( \begin{array}{cc} 
\alpha\mu M^{3} & 0 \\ 2^{-1}\mu(2\alpha-D)M^{3} & 4\alpha\mu M^{3}
\end{array} \right), \\
&E_{10}: \left( \begin{array}{cc} 
-\alpha\mu M^{3} & 0 \\ -2^{-1}\mu(2\alpha-D)M^{3} & -4\alpha\mu M^{3}
\end{array} \right)
\end{align*}
Therefore, $E_{9}$ is an unstable equilibrium and $E_{10}$ is an asymptotically stable equilibrium.

\subsection{Dynamics and connecting orbits on the Poincar\'e-Lyapunov disk}
\label{sub:SKTSS-pd7}
\begin{figure}[h]
\centering
\includegraphics[width=8cm]{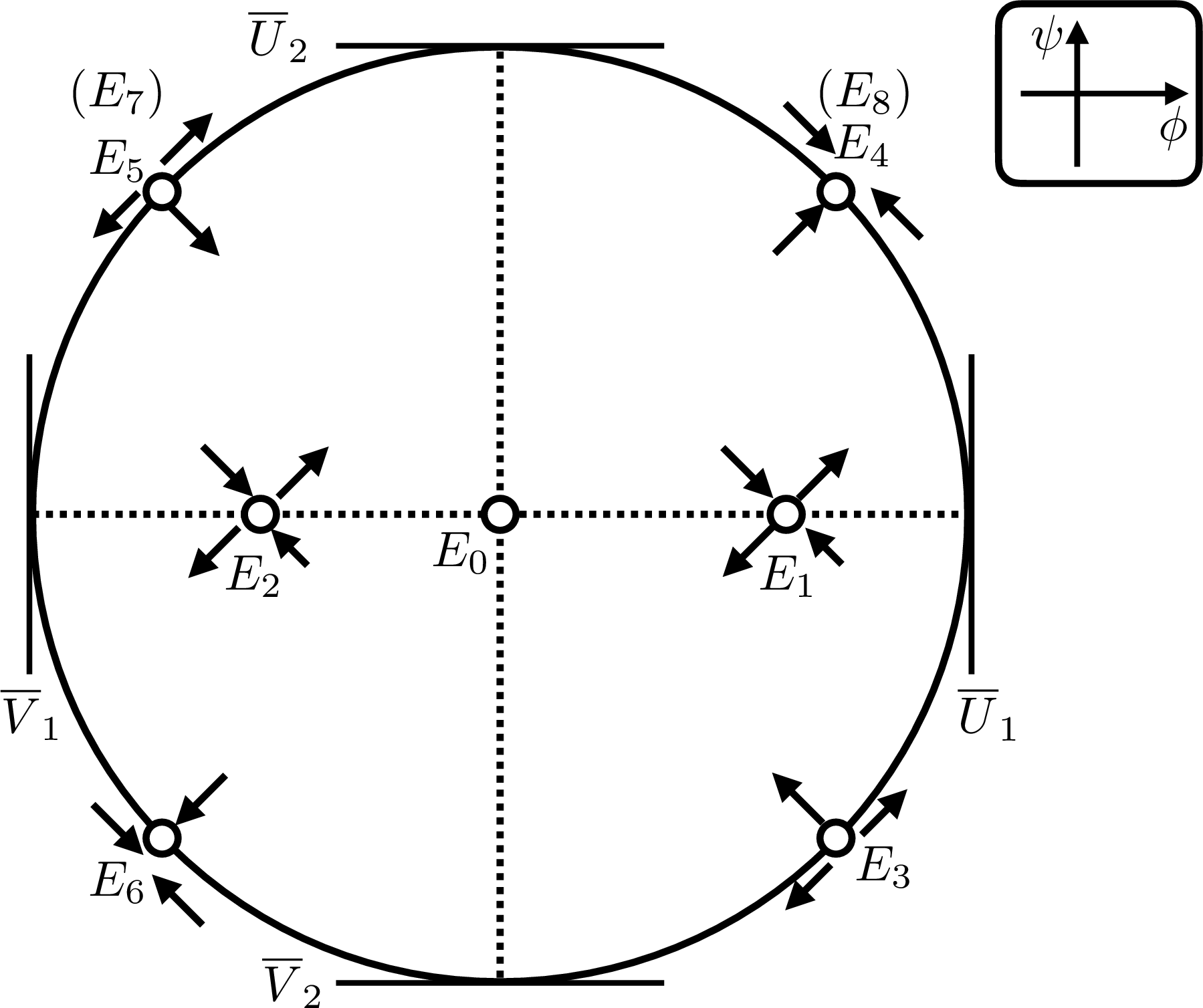}
\caption{Schematic picture of the dynamics on the Poincar\'e-Lyapunov disk.}
\label{fig:SKTSS-pd1}
\end{figure}

Combining the dynamics on the charts $\overline{U}_{j}$ and $\overline{V}_{j}$ ($j=1,2$), we obtain the whole dynamics on the Poincar\'e-Lyapunov disk, which is equivalent to the dynamics of \eqref{eq:SKTSS-pd2}.
By combining the information of finite equilibria and equilibria at infinity, we obtain Figure \ref{fig:SKTSS-pd1}.
Hereafter, the set $\Phi$ is $\Phi=\{ (u, v) \in \mathbb{R}^{2} \cup \{ \|(u, v)\|=+\infty \} \}$.

The symbols used in this section are as follows
\begin{itemize}
\item 
Let $\mathcal{W}^{s}(E_{1})$, $\mathcal{W}^{s}(E_{2})$, $\mathcal{W}^{s}(E_{4})$ and $\mathcal{W}^{s}(E_{6})$ be stable manifolds of $E_{1}$, $E_{2}$, $E_{4}$ and $E_{6}$ in the dynamical system \eqref{eq:SKTSS-pd2}.
\item 
Let $\mathcal{W}^{u}(E_{1})$, $\mathcal{W}^{u}(E_{2})$, $\mathcal{W}^{u}(E_{3})$ and $\mathcal{W}^{u}(E_{5})$ be unstable manifolds of $E_{1}$, $E_{2}$, $E_{3}$ and $E_{5}$ in the dynamical system \eqref{eq:SKTSS-pd2}.
\end{itemize}
Note that the circumference of Figure \ref{fig:SKTSS-pd1} corresponds to $\{\|(u, v)\|=+\infty\}$.

The purpose of this section is to give the existence of connected orbits in $\Phi$ and their classification results.
First, there exists a conserved quantity $H(u,v)$ in \eqref{eq:SKTSS-pd2} such that
\begin{equation}
H(u,v)= \dfrac{1}{2}\alpha \mu u^{4}+\dfrac{\mu D-2\alpha\mu}{3}u^{3}-\dfrac{\mu D}{2}u^{2}-\dfrac{1}{2}v^{2}.
\label{eq:SKTSS-pd9}
\end{equation}
That is, the function $H$ satisfies $H'(u,v)\equiv 0$ (i.e., $H$ is the first integral of \eqref{eq:SKTSS-pd2}).
This gives an orbit that, for a given constant $C$, the level curves $H(u, v)=C$ are a solution of \eqref{eq:SKTSS-pd2} in the bounded domain of $\mathbb{R}^{2}$.

Then, when $D=2\alpha$, the following two symmetries hold in \eqref{eq:SKTSS-pd2}:
\begin{equation}
(u,v) \mapsto (u, -v), \quad x\mapsto -x
\label{eq:SKTSS-pd10}
\end{equation}
and 
\begin{equation}
(u,v) \mapsto (-u, v), \quad x\mapsto -x.
\label{eq:SKTSS-pd11}
\end{equation}
Using the conserved quantities \eqref{eq:SKTSS-pd9} and the symmetries \eqref{eq:SKTSS-pd10} and \eqref{eq:SKTSS-pd11}, and the Poincar\'e-Bendixson theorem (for instance, see \cite{Wiggins} and references therein), we can conclude that there are connecting orbits and closed orbits in the case $D=2\alpha$.

That is, consider the case where an orbit starting from a point on $\mathcal{W}^{u}(E_{1})$ in $\{v<0\}$ goes to a point on $\mathcal{W}^{s}(E_{2})$, and an orbit starting from a point on $\mathcal{W}^{u}(E_{2})$ in its symmetric $\{v>0\}$ goes to a point on $\mathcal{W}^{s}(E_{1})$.
So we know that there is a periodic orbit in the neighborhood of $E_{0}$.
Also, by the nullcline on the $u$-axis and the Poincar\'e-Bendixson theorem, there are trajectories starting from an unstable manifold $\mathcal{W}^{u}(E_{3})$ at an asymptotically stable equilibrium $E_{3}$ to points on $\mathcal{W}^{s}(E_{1})$, $\mathcal{W}^{s}(E_{4})$ and $\mathcal{W}^{s}(E_{6})$. 
This implies the existence of connecting orbits between $E_{3}$ and $E_{1}$, $E_{3}$ and $E_{4}$, and $E_{3}$ and $E_{6}$.
The existence of the remaining connecting orbits is also shown in the same way by symmetry.

\begin{figure}[h]
\centering
\includegraphics[width=7cm]{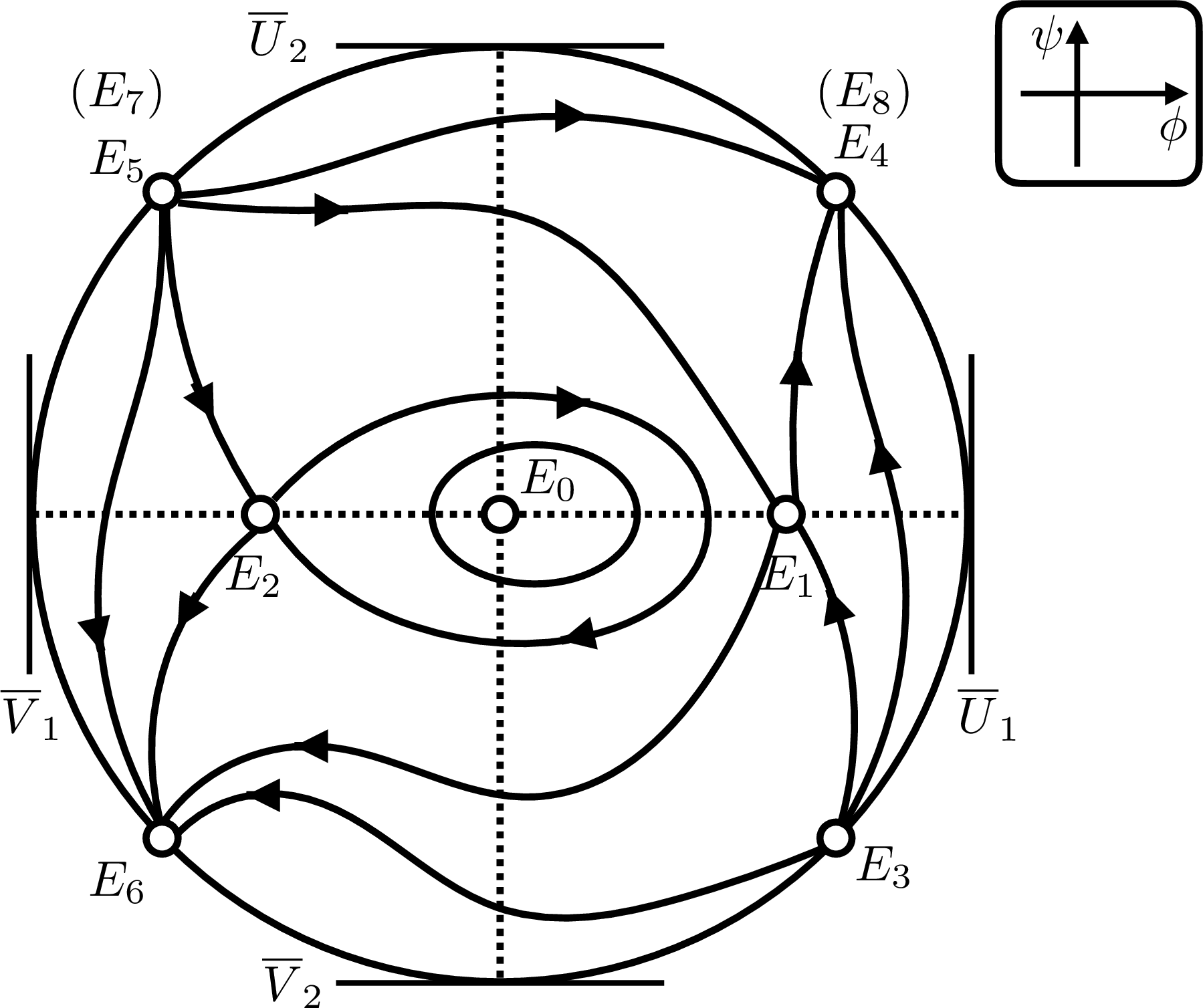}
\includegraphics[width=7cm]{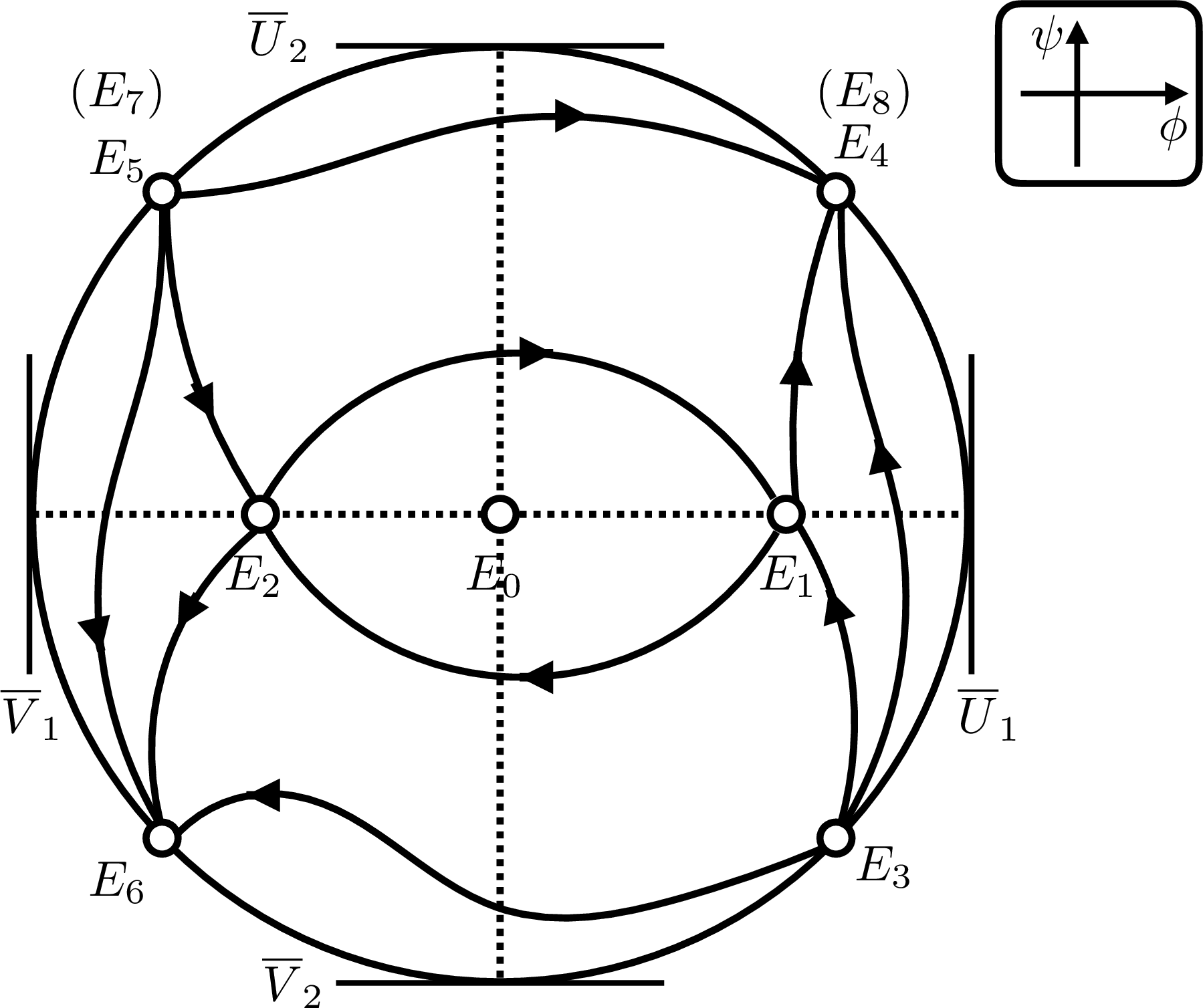}
\includegraphics[width=7cm]{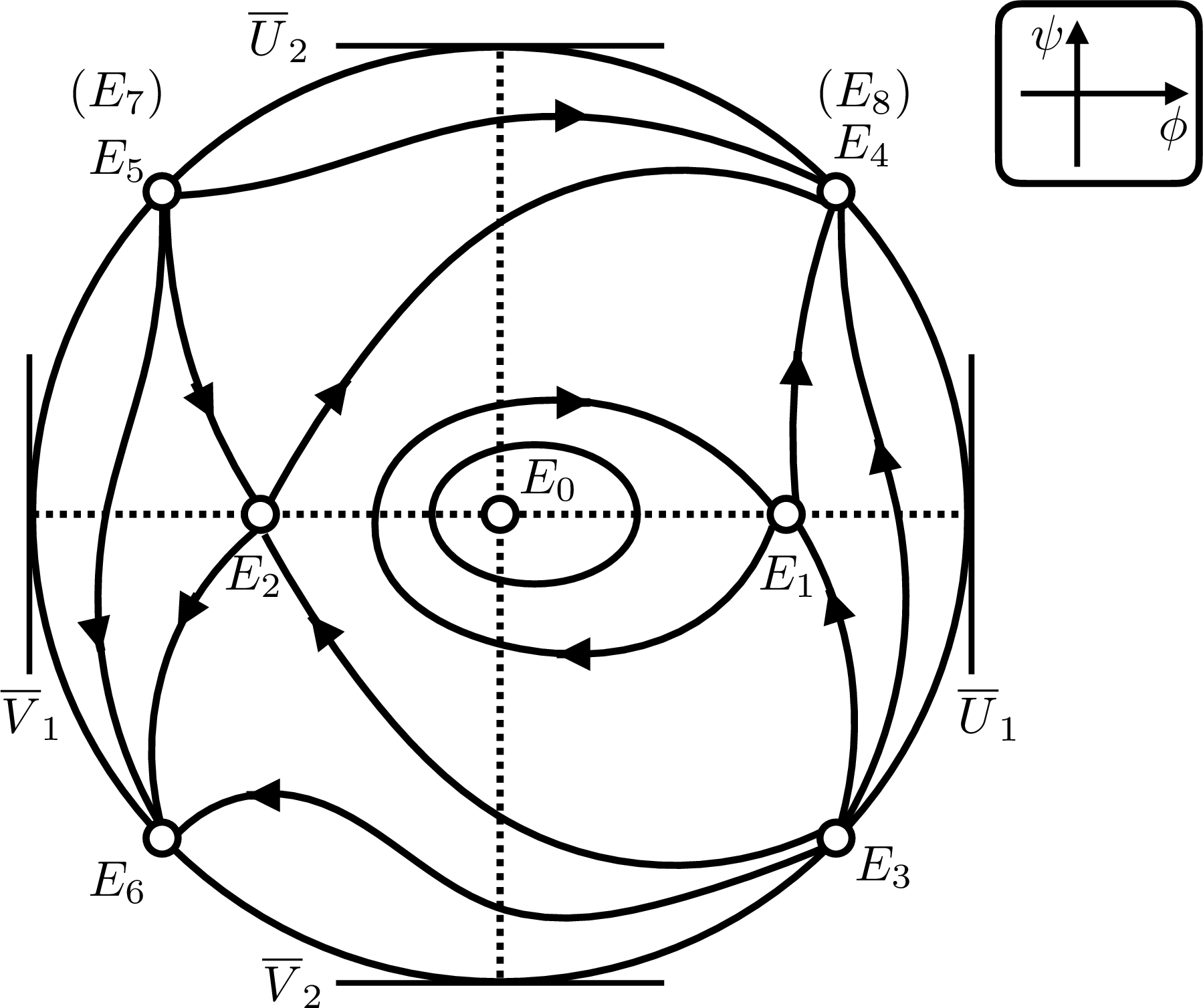}
\caption{Schematic pictures of the dynamics and connecting orbits on the Poincar\'e-Lyapunov disk for \eqref{eq:SKTSS-pd2}.
[Top left: Case $D<2\alpha$].
[Top right: Case $D=2\alpha$].
[Lower center: Case $D>2\alpha$].}
\label{fig:SKTSS-pd2}
\end{figure}

Finally, consider the case of $D>2\alpha$.
The transformation from $D>2\alpha$ to $D<2\alpha$ corresponds to the symmetry \eqref{eq:SKTSS-pd11}.
Therefore, it is sufficient to show only the case $D>2\alpha$.
When $D>2\alpha$, \eqref{eq:SKTSS-pd2} has the symmetry \eqref{eq:SKTSS-pd10}.
If a trajectory starting from a point on $A$ cannot reach a point on $B$, we can conclude that it can only go to a point on $\mathcal{W}^{s}(E_{1})$ in $\{v>0\}$ by the symmetry\eqref{eq:SKTSS-pd2}.
This means that there exists an orbit from $E_{1}$ to $E_{1}$ and a periodic orbit in the neighborhood of $E_{0}$.
Similar to the above discussion, the existence of connecting orbits between $E_{3}$ and $E_{6}$, $E_{3}$ and $E_{2}$, $E_{3}$ and $E_{1}$, and $E_{3}$ to $E_{4}$ is shown.
The symmetry \eqref{eq:SKTSS-pd10} can be used to easily show the existence of the remaining connecting orbits.
From the above, we give the existence and classification results of all connecting orbits including to infinity of \eqref{eq:SKTSS-pd2} for $D<2\alpha$, $D=2\alpha$ and $D>2\alpha$, respectively.
See Figure \ref{fig:SKTSS-pd2} for the classification results.

\begin{figure}[h]
\centering
\includegraphics[width=7cm]{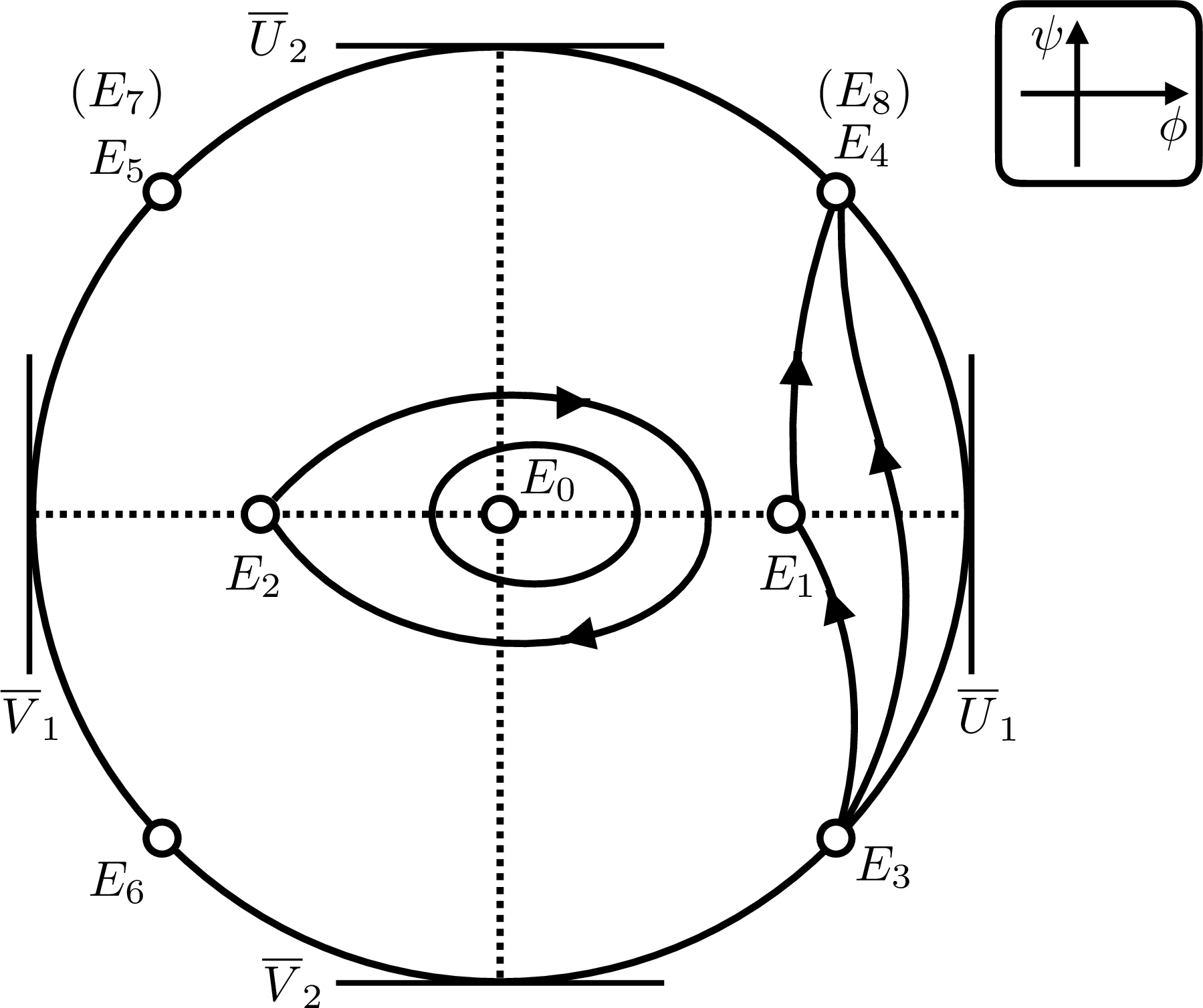}
\includegraphics[width=7cm]{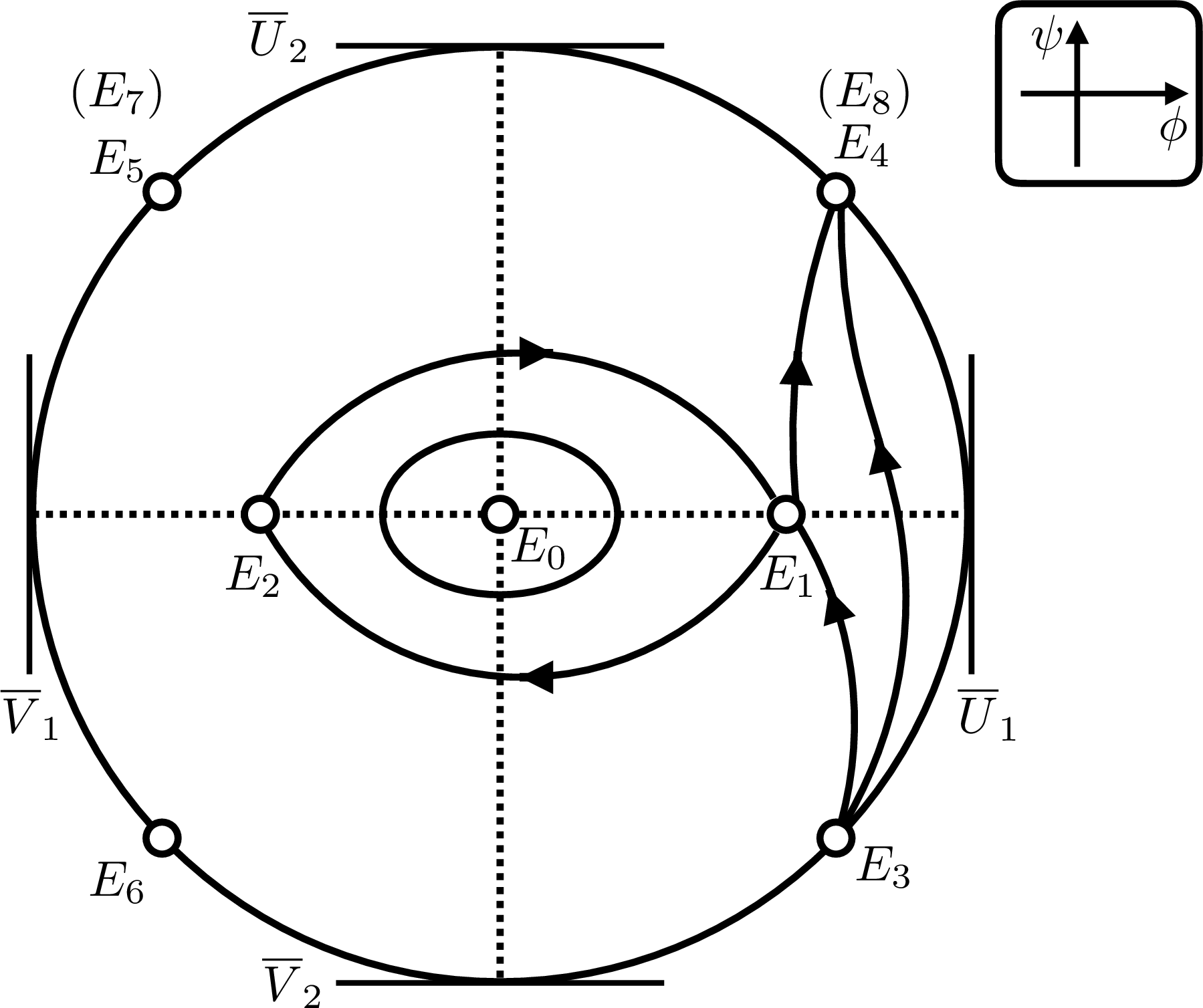}
\includegraphics[width=7cm]{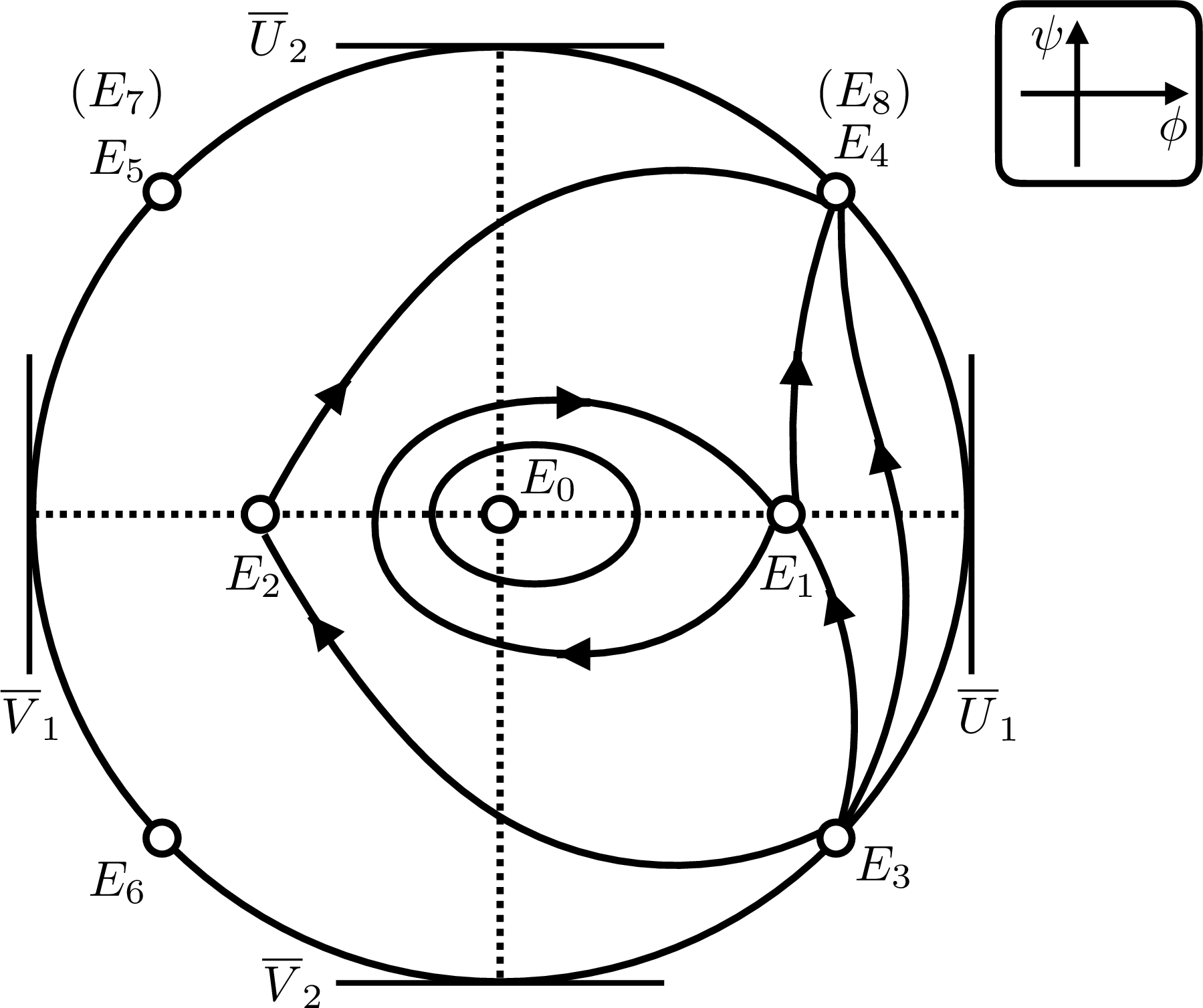}
\caption{Schematic pictures of the dynamics and connecting orbits on the Poincar\'e-Lyapunov disk for \eqref{eq:SKTSS-int8}.
[Top left: Case $D<2\alpha$].
[Top right: Case $D=2\alpha$].
[Lower center: Case $D>2\alpha$].}
\label{fig:SKTSS-pd3}
\end{figure}

The above classification of connecting orbits is for \eqref{eq:SKTSS-pd2}, and the remainder of this section will discuss the classification of connecting orbits in \eqref{eq:SKTSS-int8}.
In Figure \ref{fig:SKTSS-pd2}, we need to be careful about the treatment of the point $E_{2}$.
A note analogous to this treatment was also given in \cite{cDNPE, sDNPE, pDNPE, IM-DNPE}, but is reproduced for the reader's convenience.
When we consider the parameter $\tau$ on the disk, $E_{2}$ is the equilibrium of \eqref{eq:SKTSS-pd2}. 
However, under the parameter $x$, $E_{2}$ is a point on the line $\{D+2\alpha u=0\}$ with singularity.
We see that $du/dv$ takes the same values for \eqref{eq:SKTSS-int8} and \eqref{eq:SKTSS-pd2} except for the singularity $\{D+2\alpha u=0\}$, which implies that the trajectories for the parameter $\tau$ on $\{(u, v) \mid u\ge -(2\alpha)^{-1}D\}$ are the same as those for the parameter $x$. 
Therefore, if we focus on the stationary states in $\{u\ge -(2\alpha)^{-1}D\}$, we can see that the classification of the connecting orbits corresponding to the solutions in \eqref{eq:SKTSS-pd2} corresponds to that of the solutions in \eqref{eq:SKTSS-int8}.
The classification of the connecting orbits for $\{u\ge -(2\alpha)^{-1}D\}$ in \eqref{eq:SKTSS-int8} is shown in Figure \ref{fig:SKTSS-pd3}.

\section{Proof of theorems}
\label{sec:SKTSS-pro}
\subsection{Proof of Theorem \ref{th:SKTSS-mr1}}
\label{sub:SKTSS-pro1}
The existence of functions satisfying \eqref{eq:SKTSS-int6} in each case corresponds to the respective connecting orbits shown in Subsection ref{sub:SKTSS-pd7}.
(I) corresponds to the connecting orbit between $E_{1}$ and $E_{4}$, (II) to the connecting orbits between $E_{3}$ and $E_{4}$, (III) to the connecting orbit between $E_{3}$ and $E_{1}$, (IV) to the family of periodic orbits near $E_{0}$.

First, we derive \eqref{eq:SKTSS-mr1}.
The derivation policy is the same as for \cite{BIRD}.
The solution around $E_{4}$ are approximated as 
\[
\left(\begin{array}{cc}
\lambda_{1}(s) \\ \lambda_{2}(s)
\end{array}\right)
= C_{1}e^{-\sqrt{\alpha\mu}s}\left(\begin{array}{cc}
3\sqrt{\alpha\mu} \\ -2\alpha\mu+D\mu
\end{array}\right)
+ C_{2}e^{-4\sqrt{\alpha\mu}s}\left(\begin{array}{cc}
0 \\ 1
\end{array}\right)
+\left(\begin{array}{cc}
0 \\ \sqrt{\alpha\mu}
\end{array}\right)
\]
with constants $C_{1,2}$.
Using this equation and the some time-rescaling, we obtain
\[
\dfrac{ds}{dx}
= \dfrac{ds}{d\tau}\dfrac{d\tau}{dx}
= \lambda_{1}^{-1}\cdot \dfrac{1}{D+2\alpha u} = \dfrac{1}{D\lambda_{1}+2\alpha} \sim (2\alpha)^{-1}
\]
This yields $x(s)=2\alpha s+C_{3}$ with a constant $C_{3}$.
This means that  $s\to +\infty$ corresponds to $x\to +\infty$.
Hence, we have
\[
u(x) =\lambda_{1}^{-1}
\sim \left( 3\sqrt{\alpha\mu}C_{1}e^{-\sqrt{\alpha\mu}s}\right)^{-1}
\sim C_{4}e^{\frac{\sqrt{\alpha\mu}}{2\alpha}x}
\quad {\rm{as}} \quad x\to +\infty.
\]
The derivation of \eqref{eq:SKTSS-mr1} is completed.

Next, \eqref{eq:SKTSS-mr2} is shown.
From the subsection \ref{sub:SKTSS-pd1}, the eigenvalue of $E_{1}$ is $\omega_{\pm}$ and the target eigenvector is $(1, \omega_{+})^{T}$.
Therefore, the approximation of the solution that converges to $E_{1}$ along $\mathcal{W}^{u}(E_{1})$ as $\tau\to -\infty$ is 
\[
u(\tau) \sim 1+C_{5}e^{\omega_{+}\tau}
\quad {\rm{as}} \quad \tau \to -\infty
\]
with a constant $C_{5}$.
Here, by using the time scale transformation \eqref{eq:SKTSS-pd1}, we obtain $dx/d\tau\sim D+2\alpha$ since 
\[
\dfrac{d\tau}{dx}=\dfrac{1}{D+2\alpha u}\sim \dfrac{1}{D+2\alpha}
\quad {\rm{as}} \quad \tau\to -\infty
\]
holds.

Given that this transformation does not change the direction of time in $\{u\ge-(2\alpha)^{-1}D\}$, $\tau\to -\infty$ corresponds to $x\to -\infty$.
Hence, we obtain
\[
u(x) \sim 1+C_{5}e^{\omega_{+}\tau}
\sim 1+C_{6}e^{\frac{\omega_{+}}{2\alpha+D}x}
\quad {\rm{as}} \quad x\to -\infty.
\]
and \eqref{eq:SKTSS-mr2}.
Similar arguments also lead to \eqref{eq:SKTSS-mr3} and \eqref{eq:SKTSS-mr4}.

Then, the trajectory from $E_{3}$ to $E_{4}$ corresponding to the case (II) from the discussion in Subsection \ref{sub:SKTSS-pd7} passes through $v=0$ (i.e., the $u$-axis) only once.
Therefore, it can be seen that there exists only one $x=x_{0}$ such that $v(x)=0$.
The above completes the proof of Theorem \ref{th:SKTSS-mr1}.
\qed
\\

\subsection{Proof of Theorem \ref{th:SKTSS-mr2}}
\label{sub:SKTSS-pro2}
The existence of a function satisfying \eqref{eq:SKTSS-int6} corresponds to the ``big'' saddle homoclinic orbit from $E_{2}$ to $E_{2}$ shown in Subsection \ref{sub:SKTSS-pd7}.
The existence of $x^{*}$ is indicated by the fact that this orbit passes through the $u$-axis.
The other things to show are \eqref{eq:SKTSS-mr5} and \eqref{eq:SKTSS-mr6}.
We will show \eqref{eq:SKTSS-mr6}.
Since the eigenvalues of $E_{2}$ are $\Lambda_{\pm}$ and the eigenvector under consideration is $(1, \Lambda_{-})^{T}$, the approximation of the solution converging to $E_{2}$ along $\tau\to +\infty$ and $\mathcal{W}^{s}(E_{2})$ is 
\[
u(\tau) \sim C_{1}e^{\Lambda_{-}\tau}-\dfrac{D}{2\alpha}
\quad {\rm{as}} \quad \tau\to +\infty
\]
where $C_{1}>0$ is constant.
Then
\[
\dfrac{d\tau}{dx}=\dfrac{1}{D+2\alpha u}
\sim C_{2}e^{-\Lambda_{-}\tau}
\quad {\rm{as}} \quad \tau\to +\infty
\]
holds.
Setting
\[
x_{+}=\lim_{\tau \to +\infty}x(\tau),
\]
we then have 
\[
x_{+} =\int_{0}^{\infty} C_{3}e^{\Lambda_{-}\tau}\, d\tau <+\infty
\]
Therefore, we obtain that the asymptotic form
\[
x_{+}-x \sim C_{4}e^{\Lambda_{-}\tau}
\quad {\rm{as}} \quad \tau\to +\infty
\]
holds with a positive constant $C_{4}$.
This argument has been made in \cite{QTW, cDNPE, pDNPE, sDNPE} and is similar to those.
We then obtain
\[
u(x) 
\sim C_{1}e^{\Lambda_{-}\tau}-\dfrac{D}{2\alpha}
\sim C_{5}(x_{+}-x) -\dfrac{D}{2\alpha}
\quad {\rm{as}} \quad x\to x_{+}-0
\]
with a positive constant $C_{5}>0$.
Therefore, \eqref{eq:SKTSS-mr6} is shown.
By the same argument, \eqref{eq:SKTSS-mr5} is also shown.
This completes the proof of Theorem \ref{th:SKTSS-mr2}.
\qed
\\

\subsection{Proof of Theorem \ref{th:SKTSS-mr3}, \ref{th:SKTSS-mr4}}
\label{sub:SKTSS-pro3}
The existence of functions satisfying \eqref{eq:SKTSS-int6} corresponds to the existence of connecting orbits shown in Subsection \ref{sub:SKTSS-pd7}.
For Theorem \ref{th:SKTSS-mr3}, it corresponds to the connecting orbits from $E_{1}$ to $E_{2}$ and from $E_{2}$ to $E_{1}$.
The derivation of the asymptotic forms \eqref{eq:SKTSS-mr7} and \eqref{eq:SKTSS-mr8} is almost the same as in the above discussion.
For Theorem \ref{th:SKTSS-mr4} corresponds to the ``big'' saddle homoclinic orbit from $E_{1}$ to $E_{1}$, the connecting orbit from $E_{3}$ to $E_{2}$ and from $E_{2}$ to $E_{4}$.
The existence of the point $x_{*}$ is indicated by the fact that this orbit passes through the $u$-axis.
The above completes the proof of Theorem \ref{th:SKTSS-mr3}, \ref{th:SKTSS-mr4}.
\qed
\\

\section{Discussion}
\label{sec:SKTSS-di}
The stationary solution is one of the typical special solutions of partial differential equations and plays an important role in the study of the solution structure.
In this paper, the role of self-diffusion is discussed by considering a stationary problem for a parabolic equation in 1D space with a logistic growth term that takes into account self-diffusion, one of the nonlinear diffusions.
Understanding the role of self-diffusion leads to understanding complex nonlinearities such as the SKT cross-diffusion equations, and we believe that the classification results obtained in this paper for nonconstant stationary states are significant in understanding them.
In this paper, all dynamical systems including to infinity in $\{u\ge -(2\alpha)^{-1}D\}$ of \eqref{eq:SKTSS-int8}, which is the stationary problem of \eqref{eq:SKTSS-int1}, are studied by Poincar\'e-type compactification.
This leads to the enumeration and characterization of all existing nonconstant stationary states.
However, the stability of the obtained nonconstant stationary states is difficult to obtain by the method in this paper and requires a separate discussion.
This is an open problem.
The results obtained in this paper focus only on the structure of the equations, and we believe that they will contribute to a deeper insight into the solution structure as an evolution equation, such as a priori evaluation.

When $D>0$ and $\alpha>0$, the four types of functions given by Theorem \ref{th:SKTSS-mr1} are always present.
Furthermore, for self-diffusion coefficients $\alpha>0$, if $D<2\alpha$, there exists one function that satisfies the stationary problem for Theorem \ref{th:SKTSS-mr2}, in addition to these four types.
If $D=2\alpha$, there are two types of functions of Theorem \ref{th:SKTSS-mr3}, and if $D>2\alpha$, there are three types of functions of Theorem \ref{th:SKTSS-mr4}.
By focusing on the symmetries and conserved quantities of the two-dimensional ordinary differential equations \eqref{eq:SKTSS-pd2}, the authors also give a complete answer for the qualitative change of the non-constant stationary states, and give an explicit result that there is a relation of change boundary  $D=2\alpha$ for the self-diffusion coefficient $\alpha$ and the linear diffusion coefficient $D$.
This boundary is the loop of the heteroclinic orbit in \eqref{eq:SKTSS-int8} at $D=2\alpha$, becomes $D\neq 2\alpha$, loses symmetry \eqref{eq:SKTSS-pd11}, and changes to a ``big'' saddle homoclinic orbit.
It can be seen that the self-diffusion coefficient plays the role of a bifurcation parameter that characterizes the qualitative change in the solution structure of stationary problems for partial differential equations.
These results imply that the characteristics of the function satisfying the stationary problem change dynamically depending on the explicit large/small relationship between the linear diffusion coefficient $D$ and the self-diffusion coefficient $\alpha$.

\section*{Acknowledgments}
This work was partially supported by JSPS KAKENHI Grant Number 25K17306.


\appendix
\section*{Appendix A: Overview of the Poincar\'e-type compactification}
\label{appendix:ap1}
\setcounter{figure}{0}
\renewcommand{\thefigure}{A.\arabic{figure}}
The Poincar\'e-type compactification is one of the compactifications of the original phase space (the embedding of $\mathbb{R}^{n}$ into the unit upper hemisphere of $\mathbb{R}^{n+1}$).
In this appendix, we briefly introduce the Poincar\'e-type compactification. 
It should be noted that we refer \cite{FAL, QTW, pDNPE,IM-DNPE, Matsue1, Matsue2}.
Below, we discuss only the case $n=2$, which is necessary for the discussion in this paper. 
The argument is similar in the general dimension.

Let 
\[
\dot{u} = P(u,v),\quad
\dot{v} = Q(u,v),
\]
be equations on $\mathbb{R}^{2}$ where $\dot{~\mbox{}~}$ denotes $d/dt$.
The sphere $\mathbb{S}^{2} = \{ y \in \mathbb{R}^{3} \, |\, y_{1}^{2} + y_{2}^{2}+y_{3}^{2}=1\}$ at $\mathbb{R}^{3}$ defined by $(y_{1},y_{2},y_{3})$ is considered.
We divide the sphere into
\[
H_{+} = \{ y \in \mathbb{S}^{2}\,|\,y_{3}>0\},
\]
and
\[\mathbb{S}^{1} = \{y \in \mathbb{S}^{2}\, | \, y_{3}=0\}.\]
We also consider $\mathbb{R}^{2}$ as the plane in $\mathbb{R}^{3}$ defined by $(y_{1},y_{2},y_{3})=(u,v,1)$.

Let us consider the map given by
\[ 
f^{+}:\mathbb{R}^{2} \to H_{+}\, \cup\, \mathbb{S}^{1},
\]
where
\[f^{+}(\phi,\psi):= \pm \left( \dfrac{u}{\Delta(u,v)},\dfrac{v}{\Delta(u,v)},\dfrac{1}{\Delta(u,v)} \right) \]
with $\Delta(u,v) = \sqrt{u^{2}+v^{2}+1}$.
Geometrically, it is an embedding of $\mathbb{R}^{2}$ into $\mathbb{R}^{3}$ in the unit upper hemisphere, so that the infinity corresponds to $\mathbb{S}^{1}$. 
It is sufficient to consider the dynamics on $(u, v) \in \mathbb{R}^{2}\, \cup\, \{\|(u, v)\|=+\infty\}$ corresponding to $H_{+}\cup\mathbb{S}^{1}$, which is called Poincar\'e disk.

Then we consider six local charts on $\mathbb{S}^{2}$ given by $U_{k} = \{y \in \mathbb{S}^{2} \, | \, y_{k}>0\}$, $V_{k} = \{y \in \mathbb{S}^{2} \, | \, y_{k}<0\}$ for $k=1,2,3$.
Consider the local projection
\[ g^{+}_{k} : U_{k} \to \mathbb{R}^{2}, 
\quad g^{-}_{k} : V_{k} \to \mathbb{R}^{2} \]
defined as
\[g^{+}_{k}(y_{1},y_{2},y_{3}) = - g^{-}_{k}(y_{1},y_{2},y_{3})
 = \left(\dfrac{y_{m}}{y_{k}},\dfrac{y_{n}}{y_{k}} \right) \]
 for $m<n$ and $m,n \not= k$. 
 The projected vector fields are obtained as the vector fields on the planes
\[
\overline{U}_{k} = \{y \in \mathbb{R}^{3} \, | \, y_{k} = 1\},
\quad
\overline{V}_{k} = \{y \in \mathbb{R}^{3} \, | \, y_{k} = -1\} 
\]
 for each local chart $U_{k}$ and $V_{k}$.
 We denote by $(\lambda_{2},\lambda_{1})$ the value of $g^{\pm}_{k}(y)$ for any $k$.
 
\begin{figure}[htp]
\begin{center}
\includegraphics[width=7cm]{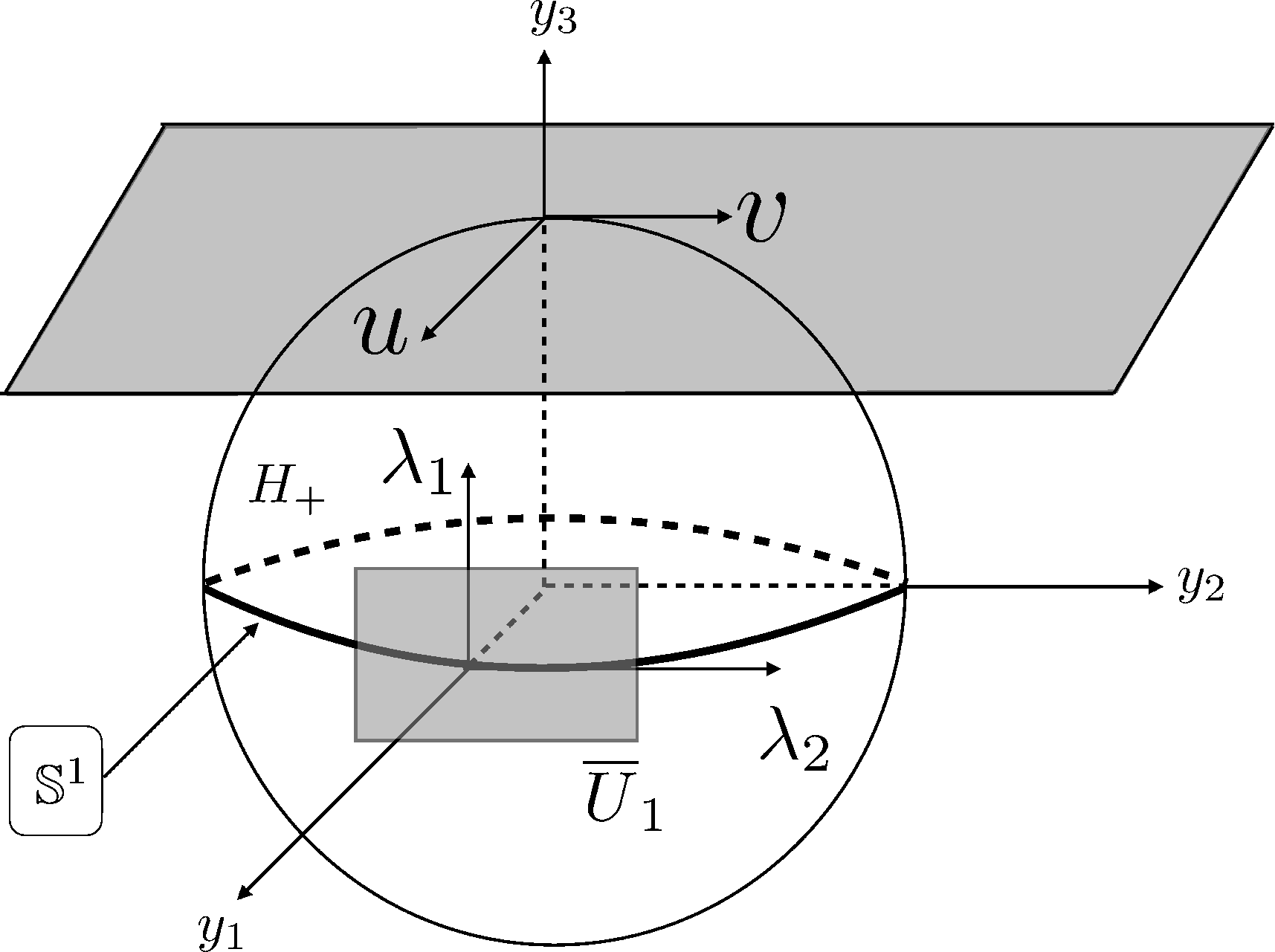}
\caption{Locations of the local chart $\overline{U}_{1}$.}
\label{fig:Poincare}
\end{center}
\end{figure}

For instance, we consider the case $k=1$.
It follows that
\[
 (g^{+}_{1} \circ f^{+})(u,v) = \left ( \dfrac{v}{u},\dfrac{1}{u}\right) = (\lambda_{2},\lambda_{1}).
 \]
Therefore, we can obtain the dynamics on the local chart $\overline{U}_{1}$ by the change of variables $u = 1/\lambda_{1}$ and $v= \lambda_{2}/\lambda_{1}$.
The locations of the Poincar\'e sphere, $(u,v)$-plane and $\overline{U}_{1}$ are expressed as Figure \ref{fig:Poincare}. 
Throughout this paper, we follow the notations used here for the Poincar\'e compactification.

Next, we consider the case that a vector field is quasi-homogeneous.
In this case, it should be noted that we choose appropriate compactifications to consider the information about dynamics at infinity.
That is, when the vector field is quasi-homogeneous, the information at infinity may not be reflected correctly in the Poincar\'e compactification.
Then, we introduce the Poincar\'e-Lyapunov compactification (the directional compactification) that is based on asymptotically quasi-homogeneous vector fields.
Then we define a class of vector fields that are quasi-homogeneous near infinity, which is determined by types and orders.
In the following, we reproduce the definitions given in \cite{Matsue1} as an aid to understanding the methods used in this paper.
See \cite{Matsue1, Matsue2} for details.
 
\begin{defn}[Case $n=2$ for \cite{Matsue1}, Definition 2.1]
 \label{def:q-h-vec}
Let $f:\mathbb{R}^{2}\to \mathbb{R}$ be a smooth function.
Let $\alpha_{1}, \alpha_{2}\ge 0$ with $(\alpha_{1}, \alpha_{2})\neq (0,0)$ be integers and $k \ge 1$. 
We say that $f$ is a quasi-homogeneous function of type $(\alpha_{1},\alpha_{2})$ and order $k$ if
\[ f(R^{\alpha_{1}}u, R^{\alpha_{2}}v)=R^{k}f(u, v),\quad \mbox{for all}\quad (u, v)\in \mathbb{R}^{2},\quad R\in \mathbb{R}. \]
Next, let 
\[
f=(f_1(u, v), f_2(u, v))
\]
be a smooth vector field. 
We say that $f=(f_1, f_2)$ is a quasi-homogeneous vector field of type $(\alpha_{1},\alpha_{2})$ and order $k+1$ if each component $f_{j}$ is a quasi-homogeneous function of type $(\alpha_{1}, \alpha_{2})$ and order $k+\alpha_{j}$ ($j=1,2$).
 \end{defn}

 \begin{defn}[Case $n=2$ for \cite{Matsue1}, Definition 2.3]
 \label{def:aqh-vec}
We say that $f$ is an asymptotically quasi-homogeneous vector field of type $(\alpha_{1}, \alpha_{2})$ and order $k+1$ at infinity if
\[
\lim_{R\to +\infty}R^{-(k+\alpha_{j})} \bigl\{ f_{j}(R^{\alpha_{1}}u, R^{\alpha_{2}}v) -R^{k+\alpha_{j}}(f_{\alpha,k})_{j}(u, v) \bigm\}=0, \quad (j=1,2)
\]
holds for any $(u,v)\in \mathbb{R}^{2}$, where $f_{\alpha,k}=((f_{\alpha,k})_{1}, (f_{\alpha,k})_{2})$ is a quasi-homogeneous vector field of type $(\alpha_{1},\alpha_{2})$ and order $k+1$.
 \end{defn}

\section*{Appendix B: Results of stationary solutions in \eqref{eq:SKTSS-int9}}
\label{appendix:ap2}
\setcounter{figure}{0}
\renewcommand{\thefigure}{B.\arabic{figure}}

In this section, we give the stationary state classification results of \eqref{eq:SKTSS-int9}.
This result follows from a similar argument from the results in this paper, which clarified the dynamical system including to infinity in \eqref{eq:SKTSS-pd2}.
The following definitions are given to classify the steady states in \eqref{eq:SKTSS-int9}.

\begin{defn}
\label{def:SKTSS-mrAp1}
Assume that $\tilde{\mu}, k$ are positive constants.
Each function $u(x)$ such that it satisfies \eqref{eq:SKTSS-int9} is defined as follows according to its behavior at $x\to +\infty$ or $x\to x_{+}-0$ ($|x_{+}|<+\infty$):
\begin{enumerate}
\item[(i)]
Define a function of type $*-1$ such that 
\[
\lim_{x \to +\infty} u(x) =1, \quad  \lim_{x \to +\infty} u_{x}(x)=0
\]
is satisfied as $x\to +\infty$.
It corresponds to the orbit in \eqref{eq:SKTSS-pd2} at $(u, v)=(1, 0)$.
\item[(ii)]
Define a function of type $*-(-k^{-1})$ such that 
\[
\lim_{x \to +\infty} u(x) =-k^{-1}, \quad  \lim_{x \to +\infty} u_{x}(x)=0
\]
is satisfied as $x\to +\infty$.
It corresponds to the orbit in \eqref{eq:SKTSS-pd2} at $(u, v)=(-D/(2\alpha), 0)=(-k^{-1}, 0)$.
\item[(iii)]
Define an unbounded function of type $*-\infty$ such that 
\[
\lim_{x \to x_{+}-0} u(x) = \lim_{x \to x_{+}-0} u_{x}(x)=+\infty
\]
is satisfied as $x\to x_{+}-0$.
It corresponds to the orbit in \eqref{eq:SKTSS-pd2} at $(u, v)=(+\infty, +\infty)$.
\item[(iv)]
Define an unbounded function of type $*- (-\infty)$ such that 
\[
\lim_{x \to x_{+}-0} u(x) = \lim_{x \to x_{+}-0} u_{x}(x)=-\infty
\]
is satisfied as $x\to x_{+}-0$.
It corresponds to the orbit in \eqref{eq:SKTSS-pd2} at $(u, v)=(-\infty, -\infty)$.
\end{enumerate}
Similarly, each function $u(x)$ such that it satisfies \eqref{eq:SKTSS-int9} is defined as follows according to its behavior at $x\to -\infty$ or $x\to x_{-}+0$ ($|x_{-}|<+\infty$):
\begin{enumerate}
\item[(v)]
Define a function of type $1-*$ such that 
\[
\lim_{x \to -\infty} u(x) =1, \quad  \lim_{x \to -\infty} u_{x}(x)=0
\]
is satisfied as $x\to -\infty$.
It corresponds to the orbit in \eqref{eq:SKTSS-pd2} at $(u, v)=(1, 0)$.
\item[(vi)]
Define a function of type $(-k^{-1})-*$ such that 
\[
\lim_{x \to -\infty} u(x) =-k^{-1}, \quad  \lim_{x \to -\infty} u_{x}(x)=0
\]
is satisfied as $x\to -\infty$.
It corresponds to the orbit in \eqref{eq:SKTSS-pd2} at $(u, v)=(-k^{-1}, 0)$.
\item[(vii)]
Define an unbounded function of type $\infty-*$ such that 
\[
\lim_{x \to x_{-}+0} u(x) = +\infty, \quad \lim_{x \to x_{-}+0} u_{x}(x)=-\infty
\]
is satisfied as $x\to x_{-}+0$.
It corresponds to the orbit in \eqref{eq:SKTSS-pd2} at $(u, v)=(+\infty, -\infty)$.
\item[(viii)]
Define an unbounded function of type $(-\infty)-*$ such that 
\[
\lim_{x \to x_{+}-0} u(x) =-\infty, \quad \lim_{x \to x_{+}-0} u_{x}(x)=+\infty
\]
is satisfied as $x\to x_{-}+0$.
It corresponds to the orbit in \eqref{eq:SKTSS-pd2} at $(u, v)=(-\infty, +\infty)$.
\end{enumerate}
\end{defn}

Under the definition of a function satisfying a stationary problem, we can give the following classification of functions satisfying stationary problems.
The proof is almost the same as in this paper.
The main difference from the main result of this paper is that the unbounded function (i.e., $u(x)$ diverges to infinity) is not $x\to \infty$ but $x\to x_{\pm}\mp0$.

\begin{thm}
\label{th:SKTSS-mrAp1}
Assume that $\tilde{\mu}, k$ are positive constants.
For any given $\tilde{\mu}, k$, the equation \eqref{eq:SKTSS-int9} has functions of types $1-\infty$, $\infty-1$, $\infty-\infty$, $-\infty-(-k^{-1})$, $(-k^{-1})-(-\infty)$ and $(-\infty)-(-\infty)$.
In particular, the asymptotic behavior of $u(x)$ for $x\to x_{+}-0$ is
\[
u(x) \sim C_{1} (x_{+}-x)^{-1} \quad {\rm{as}} \quad x\to x_{+}-0.
\]
\end{thm}

\begin{thm}
\label{th:SKTSS-mrAp2}
Assume that $\tilde{\mu}, k$ are positive constants.
If $k^{-1}<1$, in addition to the functions obtained in Theorem \ref{th:SKTSS-mrAp1}, the equation \eqref{eq:SKTSS-int9} has functions of types $(-k^{-1})-(-k^{-1})$, $(-\infty)-1$ and $1-(-\infty)$.
\end{thm}

\begin{thm}
\label{th:SKTSS-mrAp3}
Assume that $\tilde{\mu}, k$ are positive constants.
If $k=1$, in addition to the functions obtained in Theorem \ref{th:SKTSS-mrAp1}, the equation \eqref{eq:SKTSS-int9} has functions of types $(-k^{-1})-1$ and $1-(-k^{-1})$.
\end{thm}

\begin{thm}
\label{th:SKTSS-mrAp4}
Assume that $\tilde{\mu}, k$ are positive constants.
If $k^{-1}>1$, in addition to the functions obtained in Theorem \ref{th:SKTSS-mrAp1}, the equation \eqref{eq:SKTSS-int9} has functions of types $1-1$, $(-k^{-1})-\infty$ and $\infty-(-k^{-1})$.
\end{thm}

\end{document}